\documentclass[12pt,epsf]{article}
\usepackage{amssymb,amsmath,eucal,amsthm}
\usepackage{graphicx}

\def\ds{\displaystyle}

\newtheorem{proposition}{Proposition}[section]
\theoremstyle{definition}

\theoremstyle{remark}

\title{The dynamical rigid body with memory}
\author{Ion Doru Albu$^a$, Mihaela Neam\c {t}u$^b$, Dumitru Opri\c {s}$^c$}
\date{ }
\begin{document}
\maketitle

\begin{tabular}{cccccccc}
\scriptsize{$^{a}$ Department of Mathematics, Faculty of
Mathematics and Informatics,
West University of Timi\c soara,}\\
\scriptsize{Bd. V. Parvan, nr. 4, 300223, Timi\c soara, Romania, e-mail: albud@math.uvt.ro,}\\
%\scriptsize{E-mail: obundau@yahoo.es,}\\
\scriptsize{$^{b}$Department of Economic Informatics and Statistics, Faculty of Economics,West University of Timi\c soara,}\\
\scriptsize{Str. Pestalozzi, nr. 16A, 300115, Timi\c soara, Romania, e-mail:mihaela.neamtu@fse.uvt.ro,}\\
%\scriptsize{E-mail:mihaela.neamtu@fse.uvt.ro,}\\
\scriptsize{$^{c}$ Department of Mathematics, Faculty of Mathematics, West University of Timi\c soara,}\\
\scriptsize{Bd. V. Parvan, nr. 4, 300223, Timi\c soara, Romania, e-mail: opris@math.uvt.ro.}\\
%\scriptsize{E-mail: opris@math.uvt.ro}\\
\end{tabular}
\begin{abstract}In the present paper we describe the dynamics of
the revised rigid body, the dynamics of the rigid body with
distributed delays and the dynamics of the fractional rigid body.
We analyze the stationary states for given values of the rigid
body's parameters.
\end{abstract}

 \par  \noindent  {\it Key words}: rigid
body,revised rigid body,distributed delay,fractional derivative.

\section{Introduction}
\par In mechanical problems, the dynamics of the rigid body has an
important role. In many papers M. Puta analyzed the dynamics of
the rigid body with control and he obtained important results.
\par Recently, the dynamics of revised rigid body, the dynamics of
the rigid body with distributed delay and the dynamics of the
rigid body with fractional derivative have been studied. The last
two aspects represent the dynamics of the rigid body with memory.
\par Our paper studies the dynamics of the revised rigid body
obtained by a metriplectic structure which is canonically
associated. We define the dynamics of the rigid body with
distributed delays. For the Euler-Poincare dynamics of the rigid
body we analyze the linearized system in an equilibrium point. We
obtain the the existence conditions for the Hopf bifurcation with
respect to the parameter of the repartition density which defines
the distributed delay. Also, we define the fractional rigid body
dynamics using the Caputo fractional derivative.

\section{The revised differential equations for the rigid body}

The differential equations for the rigid body in $ {\mathbb R}^3$
are described by a 2-antisymmetric tensor field P and the
Hamiltonian function h given by:
\begin{equation}
P(x)=(P^{ij}(x))=\left (\begin{array}{ccc}0&x^3&-x^2\\
-x^3&0&x^1\\x^2&-x^1&0
\end{array}\right ), \quad
h=\ds\frac{1}{2}(a_1(x^1)^2+a_2(x^2)^2+a_3(x^3)^2),
\end{equation} where $(x^1,x^2,x^3)^T\in {\mathbb R}^3$, $a_i\in{\mathbb
R_+}$, $i=1,2,3$, $a_1>a_2>a_3$. These differential equations are:
\begin{equation}
\dot x(t)=P(x(t))\nabla_xh(x(t)),
\end{equation} where $\dot x(t)=(\dot x^1(t), \dot x^2(t), \dot
x^3(t))^T$ and $\nabla_xh$ is the gradient of $h$ with respect to
the canonical metric on ${\mathbb R}^3$. The differential
equations (2) have been studied by M. Puta in [8].
\par Let ${\mathbb R}^3$ be the space interpreted as the space of body
angular velocities $\Omega$ equipped with the cross product as the
Lie bracket. On this space, we consider the standard Lagrangian
 kinetic energy $L(\Omega)=\ds\frac{1}{2}I\cdot \Omega$, where
$I=diag(I_1,I_2,I_3)$ is the moment of inertial tensor, so that
the general Euler-Poincare equations become the standard rigid
body equations for a freely spinning rigid body:

\begin{equation} I\dot \Omega=(I\cdot \Omega)\times\Omega.
\end{equation}

If $M=I\cdot\Omega$ is the angular momentum, then, (3) is:

\begin{equation}\dot M=M\times \Omega.
\end{equation}

\par  If $M=(I_1x(t), I_2y(t),I_3z(t))^T$,
$\Omega=(x(t),y(t),z(t))^T$ from (4) results:
\begin{equation}
\dot x(t)=\ds\frac{I_2-I_3}{I_1}y(t)z(t),\quad \dot
y(t)=\ds\frac{I_3-I_1}{I_2}x(t)z(t),\quad  \dot
z(t)=\ds\frac{I_1-I_2}{I_3}x(t)z(t),
\end{equation}
 with $I_1>I_2>I_3$.
The differential equations (5) have been studied by M. Puta in
[9]. The revised differential equations for the rigid body given
by (2) have been studied in [5]. They are described by P, h and
the tensor fields $g=(g^{ij}(x))$, where P and h are given by (1)
and g is defined by $g(x)=(g^{ij}(x)),$
\begin{equation*}
g^{ij}(x)=\ds\frac{\partial h(x)}{\partial x^i}\ds\frac{\partial
h(x)}{\partial x^j}, i\neq j, g^{ij}(x)=-\sum_{k=1,k\neq
i}^3(\ds\frac{\partial h(x)}{\partial x^k})^2, i=1,2,3
\end{equation*} and the Casimir function of structure Poisson P is:
\begin{equation*}
c(x)=\ds\frac{1}{2}((x^1)^2+(x^2)^2+(x^3)^3).
\end{equation*}

The structure $({\mathbb R}^3,P,g,h,c)$ is called a metriplectic
manifold of second kind. The revised differential system
associated to (2) is given by:
\begin{equation}
\dot x(t)=P(x)\nabla_xh(x)+g(x)\nabla_xc(x).
\end{equation}

\begin{proposition}(i) The differential equations (2) are given
by:

\begin{equation}\dot x^1(t)\!=\!(a_2\!-\!a_3)x^2(t)x^3(t), \dot
x^2(t)\!=\!(a_3\!-\!a_1)x^1(t)x^3(t), \dot
x^3(t)\!=\!(a_1\!-\!a_2)x^1(t)x^2(t);
\end{equation}
(ii) The differential equations (6) are given by:
\begin{equation}
\begin{split}
\dot x^1(t)\!=\!(a_2\!-\!a_3)x^2(t)x^3(t)\!
+\!a_2(a_1\!-\!a_2)x^1(t)(x^2(t))^2\!+\!a_3(a_1\!-\!a_3)x^1(t)(x^3(t))^2\\
\dot x^2(t)\!=\!(a_3\!-\!a_1)x^1(t)x^3(t)\!
+\!a_3(a_2\!-\!a_3)x^2(t)(x^3(t))^2\!+\!a_1(a_2\!-\!a_1)x^2(t)(x^1(t))^2\\
\dot x^3(t)\!=\!(a_1\!-\!a_2)x^1(t)x^2(t)\!
+\!a_1(a_3\!-\!a_1)x^3(t)(x^1(t))^2\!+\!a_2(a_3\!-\!a_2)x^3(t)(x^2(t))^2;
\end{split}
\end{equation}
\end{proposition}
The equilibrium points of the system (7) are studied in [8], and
the equilibrium points of the system (8), in [5].

\section{The differential system with distributed delay for the
rigid body}

\hspace{0.6cm} We consider the space ${\mathbb R}^3$, the product
${\mathbb R}^3\times{\mathbb R}^3=\{(\tilde x,x),\tilde
x\in{\mathbb R}^3, x\in{\mathbb R}^3\}$ and the canonical
projections $\pi_i:{\mathbb R}^3\times{\mathbb R}^3\rightarrow
{\mathbb R}^3$, $i=1,2$. A vector field $X\in{\cal X}({\mathbb
R}^3\times {\mathbb R}^3)$ satisfying the condition
$X(\pi^*_1f)=0$, for any $f\in C^\infty({\mathbb R}^3)$ is given
by:
\begin{equation*}
X(\tilde x,x)=\sum_{i=1}^3X_i(\tilde
x,x)\ds\frac{\partial}{\partial x_i}.
\end{equation*}

The differential system associated to $X$ is given by:
\begin{equation*}
\dot x^i(t)=X^i(\tilde x(t),x(t)), i=1,2,3.
\end{equation*}
A differential system with distributed delay is a differential
system associates to a vector field $X\in{\cal X}({\mathbb
R}^3\times {\mathbb R}^3)$ for which $X(\pi^*_1f)=0$, for any
$f\in C^\infty({\mathbb R}^3)$ and it is given by (1), where
$\tilde x(t)$ is:
\begin{equation*}
\tilde x(t)=\int_0^k k(s)x(t-s)ds,
\end{equation*} where $k(s)$ is a density of repartition.
In what follows, we will consider the case of the following
densities of repartition:

(i) the uniform density with:
\begin{equation*}
k_\tau(s)=\left\{\begin{array}{lll} 0,& 0\leq s\leq a\\
\ds\frac{1}{\tau}, & a\leq s\leq a+\tau\\
0, &s>a+\tau\end{array},\right.
\end{equation*} where $a>0$, $\tau>0$ are given numbers;

(ii) the exponential density, with $k_\alpha(s)=\alpha e^{-\alpha
s}, \quad \alpha>0;$

(iii) the Erlang density, with $k_\alpha(s)=\alpha^2se^{-\alpha
s}, \quad \alpha>0$;

(iv) the Dirac density, with $k_\tau(s)=\delta(s-\tau), \quad
\tau>0$.

The initial condition is: $ x(s)=\varphi(s), \quad
s\in(-\infty,0],$ \\ with $\varphi:(-\infty,0]\rightarrow{\mathbb
R}^3$ a smooth map. Some systems of differential equations with
distributed delay in ${\mathbb R}^3$ were studied in [1], [2]. For
such a system, we consider relevant the geometric properties of
the vector field which defines the system, for example first
integrals (constant of the motion), Morse functions, almost
metriplectic structure, etc. The differential equations with
distributed delay for rigid body are ge\-ne\-ra\-ted by a
2-antisymmetric tensor field P on ${\mathbb R}^3\times{\mathbb
R}^3$ that satisfies the following relations:
$P(\pi^*_1f_1,\pi^*_2f_2)=0, P(\pi^*_2f_1,\pi^*_2f_2)=0,$  for all
$f_1$, $f_2$ $\in C^\infty({\mathbb R}^3)$ and $h\in
C^\infty({\mathbb R}^3\times{\mathbb R}^3)$. The differential
equation is given by:
\begin{equation}
\dot x=P(\tilde x,x)\nabla_xh(\tilde x,x).
\end{equation}

Let $P(\tilde x,x)$ be the tensor field  with the components given
by:
\begin{equation}
(P^{ij}(x,x))=\left(\begin{array}{ccc} 0 & x^3 &-\tilde
x{}^2\\
-x^3 & 0 & x^1\\
\tilde x{}^2 & -x^1 & 0\end{array}\right)
\end{equation} and
\begin{equation}
h\left(\tilde x, x\right)=a_1\tilde x{}^1x^1+a_2\tilde x{}^2
x^2+a_3\tilde x{}^3 x^3.
\end{equation}

The differential equations (9) are:
\begin{equation*}
\begin{array}{l}
\vspace{0.1cm} \dot x{}^1(t)=a_2 \tilde x^2(t)x^3(t)-a_3 \tilde
x^2(t)\tilde x^3(t),\\
\vspace{0.1cm} \dot x{}^2(t)=a_3 x^1(t)\tilde x^3(t)-a_1
\tilde x^1(t)x^3(t),\\
\dot x{}^3(t)=a_1 \tilde x^1(t)\tilde x^2(t)-a_2 x^1(t)\tilde
x^2(t).\end{array}
\end{equation*}

The differential system:
\begin{equation}
\dot x=P(\tilde x,x)\nabla_xh(\tilde x,x)+g(\tilde
x,x)\nabla_xc(\tilde x,x),
\end{equation} where the components of
$g(\tilde x,x)$ are:
\begin{equation}
g^{ij}(x)=\ds\frac{\partial h(x,\tilde{x})}{\partial
x^i}\ds\frac{\partial c(x,\tilde{x})}{\partial x^j}, i\neq j,
g^{ii}(x)=-\sum_{k=1,k\neq i}^3\ds\frac{\partial
h(x,\tilde{x})}{\partial x^k}\ds\frac{\partial
c(x,\tilde{x})}{\partial x^k}, i=1,2,3
\end{equation} is called the revised differential system with
distributed delay associated to the differential system (9).

From (10), (11), (13) and $c(\tilde
x,x)=\frac{1}{2}(x^1)^2+x^2\tilde x^2+\frac{1}{2}(x^3)^2$ results:
\begin{equation*}
(g^{ij})\!=\!\left(\begin{array}{cccc} \vspace{0.1cm}
 -a_2^2x^2\tilde x{}^2-a_3x^3\tilde
x{}^3 & a_1a_2\tilde x{}^1 x^2 & a_1 a_3\tilde x{}^1 x^3\\
\vspace{0.1cm} a_1a_2\tilde x{}^1x^2 & -a_1^2x^1\tilde
x{}^1-a_3x^3\tilde
x{}^3 & a_2a_3\tilde x{}^2x^3\\
a_1a_3\tilde x{}^1x^2 & a_2a_3\tilde x{}^2x^3 & -a_1^2x^1\tilde
x{}^1-a_2^2x^2\tilde x{}^2\end{array}\right).
\end{equation*}

The differential system (12) is given by:
\begin{equation}
\begin{array}{lll}
\vspace{0.1cm} \dot
x{}^1(t)&=&a_2\tilde x^2(t)x^3(t)-a_3\tilde x^2(t)\tilde x^3(t)+a_1a_2\tilde x^1(t)(x^2(t))^2,\\
\dot x{}^2(t)&=& a_3x^1(t)\tilde x^3(t)-a_1\tilde x^1(t)x^3(t)-a_1^2x^1(t)\tilde x^1(t)x^2(t)-a_3^2x^2(t)x^3(t)\tilde x^3(t),\\
\dot x{}^3(t)&=& a_1\tilde x^1(t)x^2(t)-a_2 x^1(t)\tilde
x^2(t)+a_2a_3\tilde x^2(t)x^2x^3.\end{array}
\end{equation}

For the Dirac distribution, system (14) was analyzed in [6]. The
other types of densities will be analyzed in our future papers.

The Euler-Poincare equation for the free rigid body with
distributed delay is defined by:
\begin{equation*}
\dot M=M\times\Omega+\alpha M\times(\tilde M\times\tilde\Omega)
\end{equation*} where $M=(I_1x(t),I_2y(t),I_3z(t))^T$,
$\Omega=(x(t),y(t),z(t))^T$, $\tilde \Omega=(\tilde x(t),\tilde
y(t), \tilde z(t))^T$, $\tilde M=I\tilde\Omega$, $I_1>0$, $I_2>0$,
$I_3>0$ and $\alpha\in{\mathbb R}$.

The equilibrium points of our system are
$\Omega_1=(\frac{m}{I_1},0,0)^T$,
$\Omega_2=(0,\frac{m}{I_2},0)^T$,
$\Omega_3=(0,0,\frac{m}{I_3})^T$, $m\in{\mathbb R}^*$.

\begin{proposition} The equilibrium point $\Omega_1$ has the
following behavior:

(i) The corresponding linear system is given by:
\begin{equation*}
\dot U(t)=AU(t)+\alpha B\tilde U(t)
\end{equation*} where $U(t)=(u^1(t),u^2(t),u^3(t))^T$ and
\begin{equation*}
A=\left(\begin{array}{ccc}0& 0& 0\\
0& 0& \frac{I_3-I_1}{I_1I_2}m\\
0& \frac{I_1-I_2}{I_1I_3}m& 0
\end{array}\right ), B=\left(\begin{array}{ccc}0& 0& 0\\
0& \frac{I_2-I_1}{I_1I_2}m^2& 0\\
0& 0& \frac{I_3-I_1}{I_1I_3}m^2
\end{array}\right );
\end{equation*}

(ii) The characteristic equation is:
\begin{equation*}
\begin{split}
&\lambda[\lambda^2-\frac{\alpha
m^2}{I_1}(\frac{I_2-I_1}{I_2}+\frac{I_3-I_1}{I_3})\lambda
k^{(1)}(\lambda)+\frac{\alpha^2m^4}{I_1^2I_2I_3}(I_2-I_1)(I_3-I_1)k^{(1)}(\lambda)^2-\\
&\frac{(I_1-I_2)(I_3-I_1)}{I_1^2I_2I_3}m^2]=0;
\end{split}
\end{equation*}

(iii) On the tangent space at $\Omega_1$ to the sphere of radius
$m^2$ the linear operator given by the linearized vector field has
the characteristic equation:

\begin{equation*}
\begin{split}
&\lambda^2-\frac{\alpha
m^2}{I_1}(\frac{I_2-I_1}{I_2}+\frac{I_3-I_1}{I_3})\lambda
k^{(1)}(\lambda)+\frac{\alpha^2m^4}{I_1^2I_2I_3}(I_2-I_1)(I_3-I_1)k^{(1)}(\lambda)^2-\\
&\frac{(I_1-I_2)(I_3-I_1)}{I_1^2I_2I_3}m^2]=0;
\end{split}
\end{equation*}

(iv) If $I_1>I_2$, $I_1>I_3$ and
$k^{(1)}(\lambda)=e^{-\tau\lambda}$, $\tau>0$, for
$0\leq\tau<\tau_c$, where
\begin{equation*}
\tau_c=\ds\frac{I_1(I_3(I_1-I_2)+I_2(I_1-I_3))}{3|\alpha|m^2(I_1-I_2)(I_1-I_3)},
\end{equation*} then the equilibrium point $\Omega_1$ is
asymptotically stable.
\end{proposition}

The analysis of the equilibrium point $\Omega_1$ for the Dirac
density is given in [1].

\section{Fractional differential systems for the rigid body.}

\hspace{0.6cm} Generally speaking, the fractional derivative,
Riemann-Liouville fractional derivative and Caputo's fractional
derivative are mostly used. In the present paper we discuss the
Caputo derivative:
\begin{equation*}
D^\alpha_tx(t)=I^{m-\alpha}(\ds\frac{d}{dt})^mx(t), \quad
\alpha>0,
\end{equation*} where $m-1<\alpha\leq m$,
$(\ds\frac{d}{dt})^m=\ds\frac{d}{dt}\circ...\circ\ds\frac{d}{dt}$,
$I^\beta$ is the $\beta$th order Riemann-Liouville integral
operator, which is expressed as follows:
\begin{equation*}
I^\beta
x(t)=\ds\frac{1}{\Gamma(\beta)}\int_0^t(t-s)^{\beta-1}x(s)ds,
\quad \beta>0.
\end{equation*}

In this paper, we suppose that $\alpha\in(0,1)$.

Examples of the fractional differential systems are: the
fractional order of Chua's system, the fractional order of
Rossler's system and the fractional Duffing oscillator. The
geometrical and mechanical interpretation of the fractional
derivative is given in [7]. The geometry of fractional osculator
bundle of higher order was made in [3] using the fractional
differential forms [4].

A fractional system of differential equations with distributed
delay in ${\mathbb R}^3$ is given by:
\begin{equation}
D^\alpha_tx(t)=X(x(t),\tilde x(t)), \quad \alpha\in(0,1),
\end{equation} where $x(t)=(x^1(t),x^2(t),x^3(t))^T\in{\mathbb
R}^3$. The linearized of (15) in the equilibrium point $x_0$ ,
$(X(x_0,x_0)=0)$ is given by the following linear fractional
differential system:
\begin{equation}
D^\alpha_tu(t)=Au(t)+B\tilde u(t),
\end{equation} where $A=(\frac{\partial X}{\partial x})|_{x=x_0}$,
$B=(\frac{\partial X}{\partial \tilde x})|_{x=x_0}$.

The characteristic equation of (16) is:
\begin{equation}
\Delta(\lambda)=det(\lambda^\alpha I-A-k^{(1)}(\lambda)B),
\end{equation} where $k^{(1)}(\lambda)=\int_0^\infty k(s)e^{-\lambda
s}ds$.

From (16) we have:
\begin{proposition} ([4])
(i) If all the roots of characteristic equation
$\Delta(\lambda)=0$ have negative real parts, then the equilibrium
point $x_0$ of (15) is asymptotically stable;

(ii) If $k(s)$ is the Dirac distribution, the characteristic
equation (17) is given by:
\begin{equation*}
\Delta(\lambda)=det(\lambda^\alpha I-A-e^{-\lambda\tau}B)=0.
\end{equation*}
If $\tau=0$, $\alpha\in(0,1)$ and all the roots of the equation
$det(\lambda I-A-B)=0$ satisfy
$|arg(\lambda)|>\frac{\alpha\pi}{2}$, then the equilibrium point
$x_0$ is asymptotically stable;

(iii) If $\alpha\in(0.5,1)$ and the equation $det(\lambda
I-A-Be^{-\lambda\tau})=0$ has no purely imaginary roots for any
$\tau>0$, then the equilibrium point $x_0$ is asymptotically
stable.
\end{proposition}

For the following delayed fractional equation (see [4])
\begin{equation}
D^\alpha_tx(t)=ax(t-\tau)
\end{equation} where $\alpha\in(0,1)$, $a\in{\mathbb R}$ and
$\tau>0$ the stability condition is:

{\it If $a<0$,
$(-a)^{\frac{1}{\alpha}}\neq\frac{1}{\tau}((2k+1)\pi-\frac{\alpha}{2}\pi)$
and $(-a)^{\frac{1}{\alpha}}\neq
-\frac{1}{\tau}((2k+1)\pi-\frac{\alpha}{2}\pi)$, $k\in{\mathbb
Z}$, then the zero solution of (18) is asymptotically stable.}

For the following delayed fractional equation (see [4])
\begin{equation}
\begin{split}
&D^\alpha_tx(t)=y(t)-k_1x(t)\\
&D^\alpha_ty(t)=-(k_1+k_2)y(t)+x(t-\tau),
\end{split}
\end{equation} where $\alpha\in(0,1)$, $k_1\geq 0$, $k_2>0$,
$\tau>0$, the stability condition is:

{\it If $k_1>0$, $k_2>\frac{1}{k_1}-k$, then the zero solution of
system (19) is asymptotically stable.}

For $f\in C^{\infty}({\mathbb R}^3)$, by $D^\alpha_{x^1}f$,
$D^\alpha_{x^2}f$, $D^\alpha_{x^3}f$ we denote the Caputo partial
derivatives defined by:
\begin{equation}
D^\alpha_{x^i}f(x)=\frac{1}{\Gamma(1-\alpha)}\int_0^{x^i}\ds\frac{\partial
f(x^1,...x^{i-1},s,x^{i+1}...,x^n)}{\partial
x^i}\frac{1}{(x^i-s)^\alpha}ds, \quad i=1,2,3
\end{equation} where $(x^i)$ are the coordinate functions on ${\mathbb
R}$, and $(\frac{\partial}{\partial x^i})$, $i=1,2,3$ is the
canonical base of the vector field on ${\mathbb R}^n$.

From (20) results:
\begin{equation*}
D^\alpha_{x^i}(x^i)^\gamma=\ds\frac{(x^i)^{\gamma-\alpha}\Gamma(1+\gamma)}{\Gamma(1+\gamma-\alpha)},
D^\alpha_{x^i}(x^j)=0, i\neq j.
\end{equation*}

Let ${\cal X}^\alpha({\mathbb R}^3)$ be the module of the
fractional vector fields generated by the operators
$\{D^\alpha_{x^i}, i=1,2,3\}$ and the module ${\cal D}({\mathbb
R}^3)$ generated by 1-forms $\{d(x^i)^\alpha, i=1,2,3\}$. The
fractional exterior derivative $d^\alpha:C^\infty({\mathbb
R}^3)\rightarrow{\cal D}({\mathbb R}^3)$ is defined by:
\begin{equation*}
d^\alpha(f)=d(x^i)^\alpha D^\alpha_{x^i}(f).
\end{equation*}

Let ${\mathop{{P}}\limits^{\alpha}}\in{\cal
X}^\alpha(\mathbb{R}^{3})\times {\cal X}^\alpha(\mathbb{R}^{3})$
be a fractional 2-skew-symmetric tensor field and $d^\alpha f$,
$d^\alpha g\in{\cal D}({\mathbb R}^3)$. The bilinear map
$[\cdot,\cdot]^\alpha:C^\infty(\mathbb{R}^{3})\times
C^\infty(\mathbb{R}^{3})\rightarrow C^\infty(\mathbb{R}^{3})$
defined by:
\begin{equation}
[f,g]^\alpha={\mathop{{P}}\limits^{\alpha}}(d^\alpha f, d^\alpha
g), \forall f,g\in C^\infty (\mathbb{R}^{3})
\end{equation} is called the fractional Leibnitz bracket.

If
${\mathop{{P}}\limits^{\alpha}}={\mathop{{P}}\limits^{\alpha}}{}^{ij}D^\alpha_{x^i}\otimes
D^\alpha_{x^j}$ then, from (21) it follows that:
\begin{equation*}
[f,g]^\alpha={\mathop{{P}}\limits^{\alpha}}{}^{ij}D^\alpha_{x^i}fD^\alpha_{x^j}g.
\end{equation*}

From the properties of the fractional Caputo, results:
\begin{equation*}
\begin{split}
&[fh,g]^\alpha=\sum\limits^\infty_{k=0}\left(\begin{array}{c}\alpha\\k\end{array}\right){\mathop{{P}}\limits^{\alpha}}{}^{ij}
(D^{\alpha-k}_{x^i}f)(D^{\alpha}_{x^j}g)\left
(\frac{\partial}{\partial x^i}\right)^kh\\
&[f,gh]^\alpha=\sum\limits^\infty_{k=0}\left(\begin{array}{c}\alpha\\k\end{array}\right){\mathop{{P}}\limits^{\alpha}}{}^{ij}
(D^{\alpha}_{x^i}f)(D^{\alpha-k}_{x^j}g)\left
(\frac{\partial}{\partial x^i}\right)^kh.
\end{split}
\end{equation*}

If ${\mathop{{P}}\limits^{\alpha}}$ is skew-symmetric we say that
$({\mathbb R}^3,[\cdot,\cdot]^\alpha)$ is a fractional almost
Poisson manifold. If $\alpha\rightarrow 1$ then we obtain the
concepts from [8].

For $h\in C^\infty({\mathbb R}^3)$, the fractional almost Poisson
dynamic system is given by:
\begin{equation}
D^\alpha_tx^i(t)=[x^i(t),h(t)]^\alpha, \textrm{where}
[x^i,h]^\alpha={\mathop{{P}}\limits^{\alpha}}{}^{ij}D^\alpha_{x^j}h.
\end{equation}

Let ${\mathop{{P}}\limits^{\alpha}}$ be a 2-skew-symmetric
fractional tensor field and ${\mathop{{g}}\limits^{\alpha}}$ a
2-symmetric fractional tensor field on ${\mathbb R}^3$. We define
the bracket $[\cdot,\cdot]^\alpha:C^\infty({\mathbb R}^3)\times
C^\infty({\mathbb R}^3)\rightarrow C^\infty({\mathbb R}^3)$ by:
\begin{equation*}
[f,h]^\alpha={\mathop{{P}}\limits^{\alpha}}(d^\alpha f, d^\alpha
g)+{\mathop{{g}}\limits^{\alpha}}(d^\alpha f, d^\alpha h), \quad
f,h\in C^\infty(\mathbb{R}^{3}).
\end{equation*}

The structure $(M, {\mathop{{P}}\limits^{\alpha}},
{\mathop{{g}}\limits^{\alpha}}, [\cdot,\cdot])$ is called
fractional almost metriplectic manifold. The fractional dynamic
system associated to $h\in C^\infty(\mathbb{R}^{3})$ is
\begin{equation}
D^\alpha_tx^i(t)=[x^i(t),h(t)]^\alpha, \textrm{where}
 [x^i,h]^\alpha={\mathop{{P}}\limits^{\alpha}}{}^{ij}D^\alpha_{x^j}h+
{\mathop{{g}}\limits^{\alpha}}{}^{ij}D^\alpha_{x^j}h.
\end{equation}

If we define the bracket
$[\cdot,(\cdot,\cdot)]^\alpha:C^\infty({\mathbb R}^3)\times
C^\infty({\mathbb R}^3)\times C^\infty({\mathbb R}^3)\rightarrow
C^\infty({\mathbb R}^3)$ by:
\begin{equation*}
[f,(h_1,h_2)]^{\alpha}={\mathop{{P}}\limits^{\alpha}}(d^\alpha f,
d^\alpha h_1)+{\mathop{{g}}\limits^{\alpha}}(d^\alpha f, d^\alpha
h_2), \forall f,h_1,h_2\in C^\infty({\mathbb R}^3),
\end{equation*} then, the fractional vector field
${\mathop{{X}}\limits^{\alpha}}{}_{h_1h_2}$ defined by:
\begin{equation*}
{\mathop{{X}}\limits^{\alpha}}{}_{h_1h_2}(f)=[f,(h_1,h_2)],
\forall f\in C^\infty({\mathbb R}^3)
\end{equation*} is called the fractional almost Leibnitz vector
field associated to the functions $h_1$, $h_2\in C^\infty({\mathbb
R}^3)$. \par The fractional almost Leibnitz dynamical system is
given by:
\begin{equation*}
D^\alpha_tx^i(t)={\mathop{{P}}\limits^{\alpha}}{}^{ij}D^\alpha_{x^j}h_1+
{\mathop{{g}}\limits^{\alpha}}{}^{ij}D^\alpha_{x^j}h_2.
\end{equation*}

Let ${\mathop{{P}}\limits^{\alpha}}=(P^{ij})$,
${\mathop{{g}}\limits^{\alpha}}=(g_{ij})$ be the fractional
2-tensor fields on ${\mathbb R}^3$ and $h\in C^\infty({\mathbb
R}^3)$ given by:
\begin{equation*}
h_1=\ds\frac{1}{\Gamma(\alpha+1)}[a_1(x^1)^{\alpha+1}+a_2(x^2)^{\alpha+1}+a_3(x^3)^{\alpha+1}].
\end{equation*}

\begin{proposition} (i) The fractional dynamic system (22) is:
\begin{equation}
D^\alpha_tx^1=(a_2-a_3)x^2x^3, D^\alpha_tx^2=(a_3-a_1)x^1x^3,
D^\alpha_tx^3=(a_1-a_2)x^1x^2;
\end{equation}

(ii) The fractional dynamic system (23) is:
\begin{equation}
\begin{split}
&D^\alpha_tx^1=(a_2-a_3)x^2x^3+a_2(a_1-a_2)x^1(x^2)^2+a_3(a_1-a_3)x^1(x^3)^2\\
&D^\alpha_tx^2=(a_3-a_1)x^1x^3+a_3(a_2-a_3)x^2(x^3)^2+a_1(a_2-a_1)x^2(x^1)^2\\
&D^\alpha_tx^3=(a_1-a_2)x^1x^2+a_1(a_3-a_1)x^3(x^1)^2+a_2(a_3-a_2)x^3(x^2)^2;
\end{split}
\end{equation}

(iii) The fractional dynamic systems (24) and (25) have the
equilibrium points $M_1(m,0,0)$, $M_2(0,m,0)$, $M_3(0,0,m)$,
$m\in{\mathbb R}^*$;

(iv) The characteristic equations for (24) are:

in $M_1(m,0,0)$: $
\lambda^\alpha(\lambda^{2\alpha}+(a_1-a_3)(a_1-a_2)m^2)=0,$

in $M_2(0,m,0)$:
$\lambda^\alpha(\lambda^{2\alpha}-(a_1-a_2)(a_2-a_3)m^2)=0,$

in $M_3(0,0,m)$:
$\lambda^\alpha(\lambda^{2\alpha}+(a_1-a_3)(a_2-a_3)m^2)=0;$ (v)
The characteristic equations in $M_1(m,0,0), \, M_2(0,m,0), \,
M_3(0,0,m)$ for (25) are:
$$\lambda^\alpha(\lambda^{2\alpha}-a_1(a_2+a_3-2a_1)m^2\lambda^\alpha+(a_1-a_3)(a_1-a_2)m^2(a_1^2m^2+1))=0,$$
$$\lambda^\alpha(\lambda^{2\alpha}-a_2(a_1+a_3-2a_2)m^2\lambda^\alpha-(a_1-a_2)(a_2-a_3)m^2(a_2^2m^2+1))=0,$$
$$\lambda^\alpha(\lambda^{2\alpha}-a_3(a_1+a_2-2a_3)m^2\lambda^\alpha+(a_1-a_3)(a_2-a_3)m^2(a_3^2m^2+1))=0;$$
\end{proposition}

The above findings allow the analysis of the equilibrium points
with respect to the parameters of the characteristic equations.

Equations (24) are called the fractional equations of the rigid
body and equations (25) are called revised fractional equations of
the rigid body. If $\alpha\rightarrow 1$, the results from
Proposition 4.2 lead to results from [8].

\par For $a_{1}:=3, \,
a_{2}:=2, \, a_{3}:=1, \, \alpha=1$ the dynamics
$(x^1(t),x^2(t),x^3(t))$ of (24) is given in figure Fig.1 and for
$\alpha=0.82$ in figure Fig.2. The numerical algorithm used is
Adam-Moulton-Bashford.

\begin{center}\begin{tabular}{ccc}
\includegraphics[width=6cm]{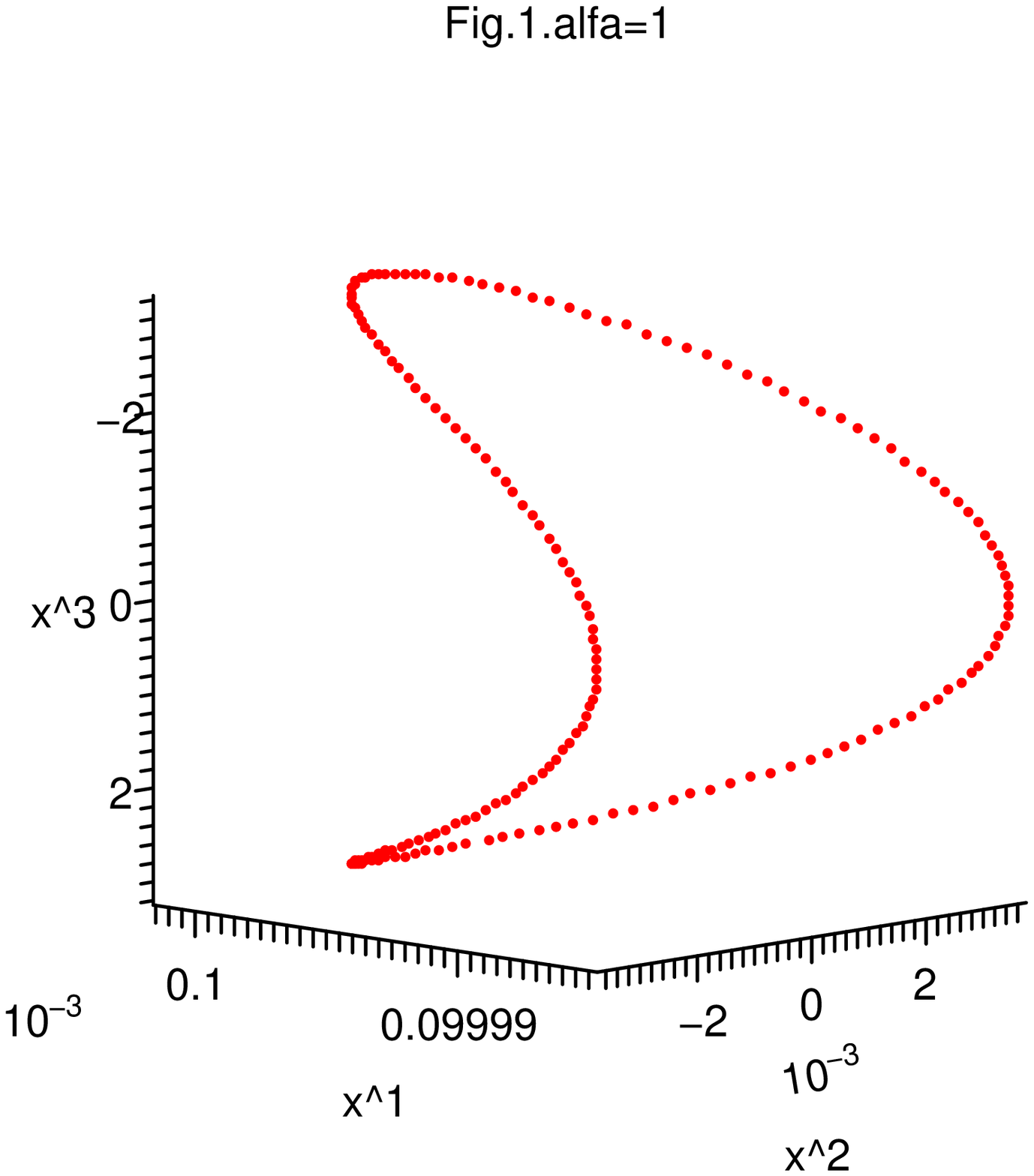} & &
\includegraphics[width=6cm]{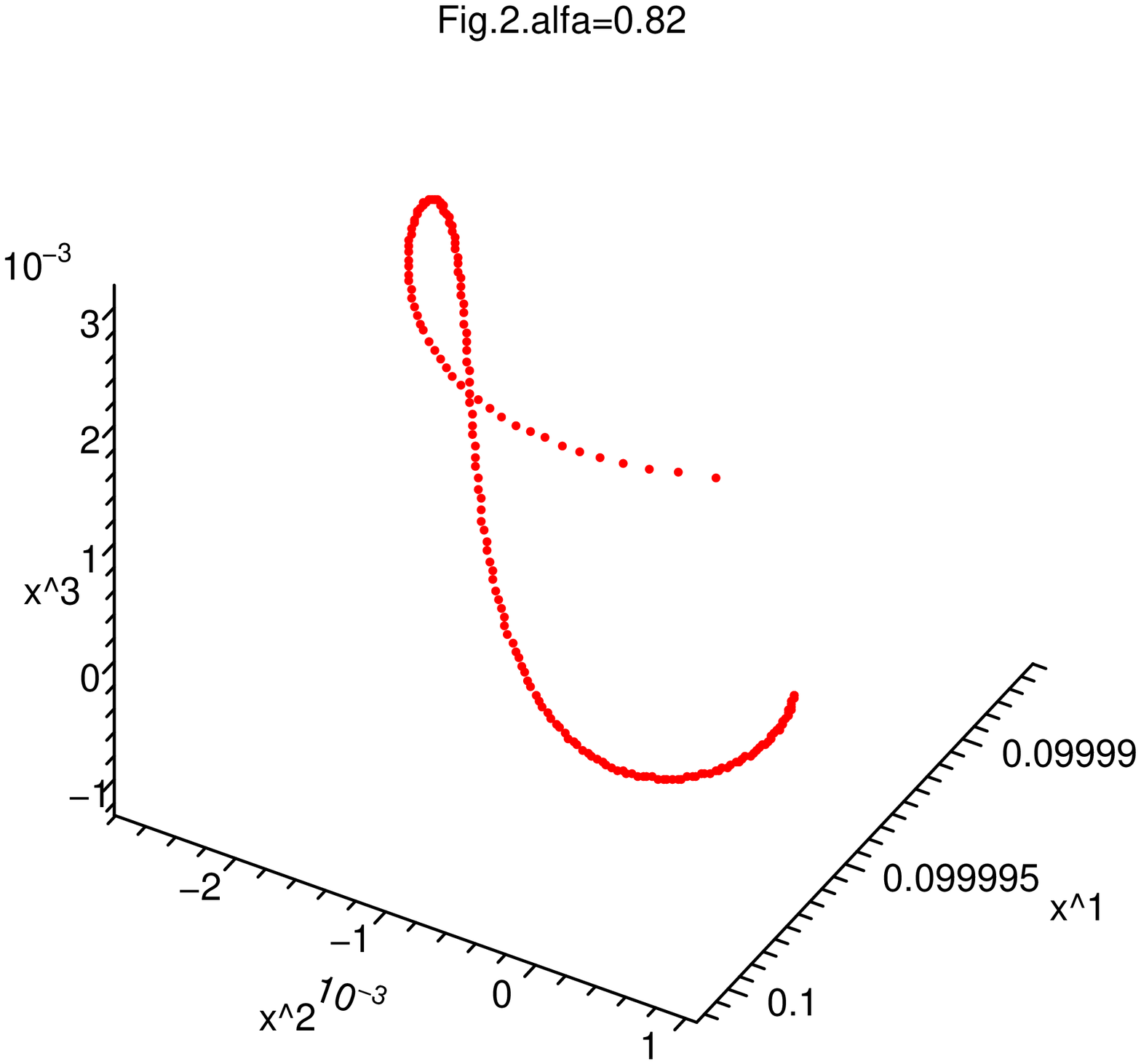}
\end{tabular}
\end{center}

\section{Conclusions}
\hspace{0.6cm} In the present paper we present the dynamics of the
rigid body with memory. The memory was described by the variables
with distributed delay and by the Caputo fractional derivative.
%The numerical simulations was realized with Maple 11.

\end{document}